\documentclass[11pt,reqno]{amsart}
\usepackage{amsthm}
\usepackage{amssymb}
\usepackage{latexsym}
\usepackage{multicol}
\usepackage{verbatim}
\usepackage{amscd}
\usepackage{hyperref}
\usepackage{amsmath, amscd}
\usepackage{color}
\usepackage[all]{xy}
\advance\textwidth by 1.2in \advance\oddsidemargin by -.6in \advance\evensidemargin by -.6in
\newtheorem*{cor}{Corollary}%[section]
\newtheorem*{lem}{Lemma}
\newtheorem*{prop}{Proposition}

\theoremstyle{definition}
\newtheorem*{defn}{Definition}
\theoremstyle{definition}
\newtheorem*{thm}{Theorem}

\newtheorem*{rem}{Remark}

\newenvironment{pf}{\proof}{\endproof}
\newcounter{cnt}
\newenvironment{enumerit}{\begin{list}{{\hfill\rm(\roman{cnt})\hfill}}{%
\settowidth{\labelwidth}{{\rm(iv)}}\leftmargin=\labelwidth%
\advance\leftmargin by \labelsep\rightmargin=0pt\usecounter{cnt}}}{\end{list}} \makeatletter
\def\mydggeometry{\makeatletter\dg@YGRID=1\dg@XGRID=20\unitlength=0.003pt\makeatother}
\makeatother \theoremstyle{remark}

\numberwithin{equation}{section}

\let\bwdg\bigwedge
\def\bigwedge{{\textstyle\bwdg}}

\begin{document}

\newcommand{\thmref}[1]{Theorem~\ref{#1}}
\newcommand{\secref}[1]{Section~\ref{#1}}
\newcommand{\lemref}[1]{Lemma~\ref{#1}}
\newcommand{\propref}[1]{Proposition~\ref{#1}}
\newcommand{\corref}[1]{Corollary~\ref{#1}}
\newcommand{\remref}[1]{Remark~\ref{#1}}
\newcommand{\defref}[1]{Definition~\ref{#1}}
\newcommand{\er}[1]{(\ref{#1})}
\newcommand{\id}{\operatorname{id}}
\newcommand{\ord}{\operatorname{\emph{ord}}}
\newcommand{\sgn}{\operatorname{sgn}}
\newcommand{\wt}{\operatorname{wt}}
\newcommand{\tensor}{\otimes}
\newcommand{\from}{\leftarrow}
\newcommand{\nc}{\newcommand}
\newcommand{\rnc}{\renewcommand}
\newcommand{\dist}{\operatorname{dist}}
\newcommand{\qbinom}[2]{\genfrac[]{0pt}0{#1}{#2}}
\nc{\cal}{\mathcal} \nc{\goth}{\mathfrak} \rnc{\bold}{\mathbf}
\renewcommand{\frak}{\mathfrak}
\newcommand{\supp}{\operatorname{supp}}
\newcommand{\Irr}{\operatorname{Irr}}
\renewcommand{\Bbb}{\mathbb}
\nc\bomega{{\mbox{\boldmath $\omega$}}} \nc\bpsi{{\mbox{\boldmath $\Psi$}}}
 \nc\balpha{{\mbox{\boldmath $\alpha$}}}
 \nc\bpi{{\mbox{\boldmath $\pi$}}}
\nc\bsigma{{\mbox{\boldmath $\sigma$}}} \nc\bcN{{\mbox{\boldmath $\cal{N}$}}} \nc\bcm{{\mbox{\boldmath $\cal{M}$}}} \nc\bLambda{{\mbox{\boldmath
$\Lambda$}}}
\newcommand{\bdd}{\operatorname{bdd}}
\newcommand{\conv}{\operatorname{conv}}
\newcommand{\lie}[1]{\mathfrak{#1}}
\makeatletter
\def\section{\def\@secnumfont{\mdseries}\@startsection{section}{1}%
  \z@{.7\linespacing\@plus\linespacing}{.5\linespacing}%
  {\normalfont\scshape\centering}}
\def\subsection{\def\@secnumfont{\bfseries}\@startsection{subsection}{2}%
  {\parindent}{.5\linespacing\@plus.7\linespacing}{-.5em}%
  {\normalfont\bfseries}}
\makeatother
\def\subl#1{\subsection{}\label{#1}}
 \nc{\Hom}{\operatorname{Hom}}
  \nc{\mode}{\operatorname{mod}}
\nc{\End}{\operatorname{End}} \nc{\wh}[1]{\widehat{#1}} \nc{\Ext}{\operatorname{Ext}} \nc{\ch}{\text{ch}} \nc{\ev}{\operatorname{ev}}
\nc{\Ob}{\operatorname{Ob}} \nc{\soc}{\operatorname{soc}} \nc{\rad}{\operatorname{rad}} \nc{\head}{\operatorname{head}}
\def\Im{\operatorname{Im}}
\def\gr{\operatorname{gr}}
\def\mult{\operatorname{mult}}
\def\Max{\operatorname{Max}}
\def\ann{\operatorname{Ann}}
\def\sym{\operatorname{sym}}
\def\Res{\operatorname{\br^\lambda_A}}
\def\und{\underline}
\def\Lietg{$A_k(\lie{g})(\bsigma,r)$}

 \nc{\Cal}{\cal} \nc{\Xp}[1]{X^+(#1)} \nc{\Xm}[1]{X^-(#1)}
\nc{\on}{\operatorname} \nc{\Z}{{\bold Z}} \nc{\J}{{\cal J}} \nc{\C}{{\bold C}} \nc{\Q}{{\bold Q}}
\renewcommand{\P}{{\cal P}}
\nc{\N}{{\Bbb N}} \nc\boa{\bold a} \nc\bob{\bold b} \nc\boc{\bold c} \nc\bod{\bold d} \nc\boe{\bold e} \nc\bof{\bold f} \nc\bog{\bold g}
\nc\boh{\bold h} \nc\boi{\bold i} \nc\boj{\bold j} \nc\bok{\bold k} \nc\bol{\bold l} \nc\bom{\bold m} \nc\bon{\bold n} \nc\boo{\bold o}
\nc\bop{\bold p} \nc\boq{\bold q} \nc\bor{\bold r} \nc\bos{\bold s} \nc\boT{\bold t} \nc\boF{\bold F} \nc\bou{\bold u} \nc\bov{\bold v}
\nc\bow{\bold w} \nc\boz{\bold z} \nc\boy{\bold y} \nc\ba{\bold A} \nc\bb{\bold B} \nc\bc{\bold C} \nc\bd{\bold D} \nc\be{\bold E} \nc\bg{\bold
G} \nc\bh{\bold H} \nc\bi{\bold I} \nc\bj{\bold J} \nc\bk{\bold K} \nc\bl{\bold L} \nc\bm{\bold M} \nc\bn{\bold N} \nc\bo{\bold O} \nc\bp{\bold
P} \nc\bq{\bold Q} \nc\br{\bold R} \nc\bs{\bold S} \nc\bt{\bold T} \nc\bu{\bold U} \nc\bv{\bold V} \nc\bw{\bold W} \nc\bz{\bold Z} \nc\bx{\bold
x} \nc\KR{\bold{KR}} \nc\rk{\bold{rk}} \nc\het{\text{ht }}
\nc\loc{\rm{loc }}

\nc\toa{\tilde a} \nc\tob{\tilde b} \nc\toc{\tilde c} \nc\tod{\tilde d} \nc\toe{\tilde e} \nc\tof{\tilde f} \nc\tog{\tilde g} \nc\toh{\tilde h}
\nc\toi{\tilde i} \nc\toj{\tilde j} \nc\tok{\tilde k} \nc\tol{\tilde l} \nc\tom{\tilde m} \nc\ton{\tilde n} \nc\too{\tilde o} \nc\toq{\tilde q}
\nc\tor{\tilde r} \nc\tos{\tilde s} \nc\toT{\tilde t} \nc\tou{\tilde u} \nc\tov{\tilde v} \nc\tow{\tilde w} \nc\toz{\tilde z}
\title{Tilting modules  for the current algebra \\  of a simple Lie algebra}
\author{Matthew Bennett and Vyjayanthi Chari}
\address{Department of Mathematics, University of California, Riverside, CA 92521,
 ]U.S.A.}\email{mbenn002@gmail.com, vyjayanthi.chari@ucr.edu}\thanks{V.C. was partially supported by DMS-0901253}
\maketitle
\begin{abstract} The category of level zero  representations of current and affine Lie algebras shares many of the properties of other well--known categories which appear in Lie theory and in algebraic groups in characteristic $p$ and  in this paper we explore further similarities. The role of the standard and co--standard module is  played by the finite--dimensional local Weyl module and the dual of the infinite--dimensional global Weyl module respectively.  We  define the canonical  filtration of a  graded module for the   current algebra. In the case when $\lie g$ is of type $\lie{sl}_{n+1}$ we show that the  well--known necessary and sufficient homological condition for a  canonical filtration to be a good (or a $\nabla$--filtration) also holds in our situation. Finally, we construct the indecomposable  tilting modules in our category and show that any tilting module is isomorphic to a direct sum of indecomposables.
\end{abstract}
\section*{Introduction}  The study of the representation theory of current algebras was largely motivated by its relationship to the representation theory of affine and quantum affine algebras associated to a simple Lie algebra $\lie g$. However,  it is also now of independent interest since it yields connections with problems arising in mathematical physics, for instance  the $X=M$ conjectures, see \cite{AK}, \cite{deFK}, \cite{Naoi}. These connections arise from the fact that the current algebra is graded by the non--negative integers and that studying graded modules and their characters give rise to interesting combinatorics. The work of \cite{KodNaoi} for instance, also relates certain graded characters to the Poincare polynomials of quiver varieties.

The current Lie algebra is just the Lie algebra of polynomial maps from $\bc\to \lie g$ and can be identified with the space $\lie g\otimes \bc[t]$ with the obvious commutator. The Lie algebra and its universal enveloping algebra inherit a  grading coming from the natural grading on $\bc[t]$. One is interested in the category $\cal I$ of $\bz$--graded modules of $\lie g[t]$ with the restriction that the graded pieces are finite--dimensional. The simple objects in the category are just the graded shifts of the irreducible modules for $\lie g$ and so are parametrized by  a set $\Lambda$ consisting of pairs $(\lambda,r)$, where $\lambda$ is a dominant integral weight and $r$ is an integer. However, the main interest of this category is that it has reducible but indecomposable objects. Many of these objects are either defined in a way similar to, or play a role  which is  analogous to well--known constructions in Lie theory, say in the BGG category $\cal O$ associated to a simple Lie algebra or to representations of algebraic groups in characteristic $p$. Our work has some similarity with  \cite{Mazor} although our set up is quite different. In particular the  grade zero piece of the algebra $\bu(\lie g[t])$ is infinite--dimensional.

The category $\cal I$ contains the projective cover and the  injective envelope of a  simple object. Moreover, if we define a suitable partial order on $\Lambda$, then we can define the appropriate analog of the standard and costandard objects in $\cal I$. An interesting feature in our case is that the standard object $\Delta(\lambda,r)$ is a finite--dimensional module called the local Weyl module which has been extensively studied (see \cite{CL1}, \cite{FoL}, \cite{Naoi2}, for instance). The co-standard object $\nabla(\lambda,r)$ however is infinite--dimensional and is the (appropriately defined) dual of the global Weyl module. Both modules lie in a nice subcategory of $\cal I$ which we call $\cal I_{\bdd}$. It is the full subcategory consisting of objects whose weights are in a finite union of cones (as in $\cal O$) and whose grades are bounded above.

The main goal of this paper is to construct another family of non--isomorphic modules indexed by $\Lambda$ and which are in $\cal I_{\bdd}$. These modules are denoted by $T(\lambda,r)$ and have an infinite filtration in which the successive quotients are of the form $\Delta(\mu,s)$ for $(\mu,s)\in\Lambda$. The filtration multiplicity of any given $\Delta(\mu,s)$ is finite. We also show that these modules satisfy a nice homological property,
namely that $$\Ext^1_{\cal I}(\Delta(\mu,s), T(\lambda,r))=0,\ \ (\mu,s), (\lambda,r)\in\Lambda.$$

In the case of algebraic groups (for instance, see \cite{Donkin},\cite{Mathieu}) it is shown that the preceding  condition is equivalent to the module having a filtration by $\nabla(\mu,s)$ and the module $T(\lambda,r)$ is then called tilting. A crucial tool in that situation to proving this equivalence is to show that every module can be embedded into a module admitting a $\nabla$--filtration.

In our case, we first have to modify slightly the definition of the $\nabla$--filtration, but the more serious problem is to show that any object embeds into one which has a $\nabla$--filtration. If we restrict our attention to $\cal I_{\bdd}$ then we are able to prove that any $M$ embeds into an injective object of $\cal I_{\bdd}$. We show that if these  injective objects admit a $\nabla$--filtration, then the modules $T(\lambda,r)$ are tilting, and are all the indecomposable tilting objects in $\cal I_{\bdd}$. Moreover, we also prove that any tilting module in $\cal I_{\bdd}$ is isomorphic to a direct sum of indecomposable tilting modules. In the case when $\lie g$ is of type $\lie{sl}_{n+1}$ (see \cite{BCM} for the $n=1$ case and \cite{BBCKL} for general $n$), it is shown that the injective envelopes of simple objects do have $\nabla$--filtrations. In fact, it is also shown in those papers that the injective envelope a simple object in  $\cal I_{\bdd}$ (which is usually smaller) also  has a $\nabla$--filtration.  As a consequence, one sees that for $\lie{sl}_{n+1}$ the modules $T(\lambda,r)$ are indeed tilting modules.  There are obviously a number of interesting questions one could ask about these modules which we will pursue elsewhere.

\section{Preliminaries}

\subsection{} Throughout this paper we denote by   $\bc$  the field of complex numbers and $\bz$ (resp. $\bz_+$) the set of integers (resp.  nonnegative   integers). For any Lie algebra $\lie a$, we denote by $\bu(\lie a)$ the universal enveloping algebra of $\lie a$.  Let $t$ be an indeterminate and let $\lie a[t]=\lie a\otimes \bc[t]$ be  the Lie algebra with commutator given by, $$[a\otimes f, b\otimes g]=[a,b]\otimes fg,\ \ a, b\in\lie a, \ \ f,g,\in\bc[t].$$  We identify $\lie a$ with the Lie subalgebra $\lie a\otimes 1$ of $\lie a[t]$. The Lie algebra $\lie a[t]$ has a natural $\bz$--grading given by the powers of $t$ and this also induces a $\bz$--grading on $\bu(\lie a[t])$, and  $$\bu(\lie a[t])[s]=0,\ \ \  s<0,\ \ \qquad  \bu(\lie a[t])[0]=\bu(\lie a).$$ The  graded pieces are  $\lie a$--module under left and right multiplication by elements of $\lie a$ and hence also under the adjoint action of $\lie a$. In particular, if $\dim\lie a<\infty$, then $\bu(\lie a[t])[r]$ is a free module for $\lie a$ (via left or right multiplication) of finite rank.

\subsection{}
From now on, $\lie g$  denotes a finite--dimensional complex simple Lie algebra of rank $n$ and $\lie h$   a fixed Cartan subalgebra of $\lie g$. Let   $I=\{1,\cdots ,n\}$  and fix a set  $\{\alpha_i: i\in I\}$ of simple roots of $\lie g$ with respect to $\lie h$ and a set $\{\omega_i: i\in I\}$   of fundamental weights.  Let $Q$ (resp. $Q^+$) be  the integer span (resp. the nonnegative integer span) of $\{\alpha_i: i\in I\}$ and similarly define  $P$ (resp. $P^+$) to be the $\bz$ (resp. $\bz_+$) span of  $\{\omega_i: i\in I\}$. Let $\{x_i^\pm, h_i: i\in I\}$ be a set of Chevalley generators of $\lie g$ and let $\lie n^\pm$ be the Lie subalgebra of $\lie g$ generated by the elements $x_i^\pm$, $i\in I$. We have, $$\lie g\ =\ \lie n^-\oplus\lie h\oplus\lie n^+,\ \ \qquad \bu(\lie g)=\bu(\lie n^-)\otimes\bu(\lie h)\otimes\bu(\lie n^+). $$

 Given $\lambda,\mu\in\lie h^*$, we say that  $\lambda\le\mu$ iff $ \lambda-\mu\in Q^+$.
Let $W$ be the Weyl group of $\lie g$ and let $w_0\in W$ be the longest element of $W$. Given $\lambda\in P^+$, let $\conv W\lambda\subset\lie h^*$ be the convex hull of the set $W\lambda$.

\subsection{}\label{locfin}
For any $\lie g$-module $M$ and $\mu\in\lie h^*$, set
\[M_\mu=\{m\in M\ :\ hm=\mu(h)m,\quad h\in\lie h\}.\] We say  $M$ is a \textit{weight module}
for $\lie g$ if \[M=\bigoplus_{\mu\in\lie h^\ast} M_\mu,\] and  we set $\wt(M)=\left\{\mu\in\lie h^\ast\ :\ M_\mu\ne 0\right\}.$
Any finite--dimensional $\lie g$--module is  a weight module.
   It is well-known that the set of isomorphism classes of irreducible finite-dimensional $\lie g$-modules is in bijective correspondence with $P^+$. For  $\lambda\in P^+$ we denote by $V(\lambda)$ a representative of the corresponding isomorphism class which is generated by  a vector $v_\lambda$ with defining relations
  \[ \lie n^+v_\lambda=0,\qquad hv_\lambda=\lambda(h)v_\lambda,\qquad (x_{{i}}^-)^{\lambda(h_{i})+1}v_{\lambda}=0,\ \ \  h\in\lie h,\ \ i\in I.\]
 and recall that $\wt V(\lambda)\subset\conv W\lambda$. The module $V(0)$ is the trivial module for $\lie g$ and we shall write it as $\bc$.
The character of $M$ is the  element of the integral group ring $\bz[P]$ defined by,
\[ {\rm ch}_{\lie g} M=\sum_{\mu\in P}\dim_\bc M_\mu e(\mu),\] where $e(\mu)\in\bz[P]$ is the generator of the group ring corresponding to $\mu$.
The  set $\{\ch_{\lie g} V(\mu): \mu\in P^+\}$ is a  linearly independent subset of $\bz[P]$.

We say that $M$ is a \textit{locally finite-dimensional} $\lie g$--module if it is a direct sum of finite--dimensional $\lie g$--modules,  in which case $M$ is necessarily a weight module.
Using Weyl's theorem one knows that  a locally finite-dimensional $\lie g$-module $M$ is isomorphic to a direct sum of modules of the form $V(\lambda)$, $\lambda\in P^+$ and  hence  $\wt M\subset P$. Set \begin{equation}\label{ndecomp} M^{\lie n^+}=\{m\in M: \lie n^+ m=0\}\ \qquad M^{\lie n^+}_\lambda=M^{\lie n^+}\cap M_\lambda\cong\Hom_{\lie g}(V(\lambda), M). \end{equation} \iffalse We have $$ M_\lambda=(\lie n^-M\cap M_\lambda)\oplus M_\lambda^{\lie n^+},\qquad
M=\bigoplus_{\lambda\in P^+}\bu(\lie g)M_\lambda^{\lie n^+},$$ or equivalently, the isotypical component of $M$ corresponding to $\lambda\in P^+$ is the $\lie g$--submodule generated by $M_\lambda^{\lie n^+}$.
\fi

\subsection{}

Let $\cal I$ be the  category whose objects are graded $\lie g[t]$-modules $V$ with finite-dimensional graded components and where the morphisms are  maps  of graded $\lie g[t]$-modules.  Thus an object $V$ of $\cal I$, is a $\bz$--graded vector space  $V =\oplus_{s\in\bz}V[s]$ which admits a left action of $\lie g[t]$ satisfying
$$(\lie g\otimes t^r)V[s]\subset V[s+r],\qquad s,r\in\bz.$$
A morphism between  two objects $V$, $W$ of $\cal I$ is a degree zero  map of graded $\lie g[t]$--modules.  Clearly $\cal I$ is closed under taking submodules, quotients and finite direct sums. For any $r\in\bz$ we let $\tau_r$ be the grade shifting operator.

If $V\in\Ob\cal I$ and $\mu\in P^+$, then 
$$V_\mu^{\lie n^+}= \bigoplus_{r\in\bz} V[r]_\mu^{\lie n^+},\ \qquad V[r]_\mu^{\lie n^+}= V_\mu^{\lie n^+}\cap V[r].$$
The \textit{graded character}  of $V\in\Ob\mathcal I$ is an element of the space of power series $\bz[P][[u,u^{-1}]]$, given by
\[{\rm ch}_{\rm{gr}}V:=\sum_{r\in\bz} {\rm ch}_\lie{g}(V[r])u^{r},\] 
where we observe that for all $r\in\bz$ the subspace $V[r]$ is a $\lie g$--module. Given $V\in\Ob\cal I$, the restricted dual is 
$$V^*=\bigoplus_{r\in\bz} V^*[r],\qquad  V^*[r]=V[-r]^*.$$  
Then $V^*\in\Ob\cal I$ with the usual action:
$$(xt^s)v^*(w)=-v^*(xt^s w),$$ 
and $(V^*)^*\cong V$ as objects of $\cal I$. Note that if $V\in\Ob{\cal I}$, then 
\[{\rm ch}_{\rm{gr}}V^*:=\sum_{r\in\bz} {\rm ch}_\lie{g}(V[r]^*)u^{-r}.\]

\section{The main result}

\subsection{} Let $\cal I_{\bdd}$ be the full subcategory of $\cal I$ consisting of objects $M$ satisfying the following two conditions:
\begin{enumerit}
\item there exists $\mu_1,\cdots, \mu_s\in P^+$ (depending on $M$) such that $$\wt M\subset \bigcup_{\ell=1}^s\conv W\mu_\ell,$$
\item there exists $r\in\bz$ (depending on $M$) such that $M[\ell]=0$ if $\ell\ge r$.
\end{enumerit} Notice that $\cal I_{\bdd}$ is not closed under taking duals.
We now define  three natural families of objects of $\cal I_{\bdd}$ which are all indexed by $P^+\times\bz$.

\subsection{} Let $\ev_0: \lie g[t]\to\lie g$ be the homomorphism of Lie algebras which maps $x\otimes f\mapsto f(0)x$. The kernel of this map is a graded ideal in $\lie g[t]$ and hence any $\lie g$--module $V$ can be regarded
as a graded $\lie g[t]$--module by pulling back through $\ev_0$ and $\ev_0 V\in\Ob\cal I$ if $\dim V < \infty$. The pull back of $V(\lambda)$ is denoted $V(\lambda,0)$ and we set $\tau_rV(\lambda,0)=V(\lambda,r)$ and we let $v_{\lambda,r}\in V(\lambda,r)$ be the element corresponding to $v_\lambda$.
For any module $M$  denote by $\soc M$  the maximal semisimple submodule of $M$. The next proposition gives an explanation for restricting our study to $\cal I_{\bdd}$.
\begin{prop}\label{simpleclass}
 \begin{enumerit}
 \item[(i)] Any irreducible  object  in $\cal I$ (or $\cal I_{\bdd}$) is isomorphic to $V(\mu,r)$ for a unique element
  $(\mu,r)\in P^+\times \bz$. Moreover $$V(\mu,r)^*\cong V(-w_0\mu,-r).$$
  \item[(ii)] Let $M\in\Ob\cal I_{\bdd}$ be non-trivial. Then $\soc M\ne 0$ and we have $$\soc M \cong \bigoplus_{(\lambda,r)\in P^+\times\bz} V(\lambda,r)^{m(\lambda, r)},\ \ m(\lambda,r)=\dim\Hom_{\cal I}(V(\lambda, r), M).$$\end{enumerit}
 \end{prop}
 \begin{pf} Part (i) is straightforward and a proof can be found in  \cite[Proposition 1.3]{CG1}.
 For (ii), choose $s\in\bz$ such that $M[s]\ne 0$ and $M[\ell]=0$ for all $\ell>s$.  Since $M[s]$ is a finite--dimensional $\lie g$--module, there exists $\mu\in P^+$ such that $\Hom_{\lie g}(V(\mu), M[s])\ne 0$. Since$$(\lie g\otimes t\bc[t])M[s]=0,$$ it follows that $\Hom_{\lie g[t]}(V(\mu,s), M)\ne 0$ proving that $\soc M\ne 0$. The rest of (ii) is now immediate.
\end{pf}

 \subsection{}\label{delta}  The next family we need are the local Weyl modules which were originally defined in \cite{CPweyl}. For the purposes of this paper, we shall denote them as $\Delta(\lambda,r)$,  $(\lambda,r)\in P^+\times\bz$. Thus, $\Delta(\lambda, r)$ is generated as a $\lie g[t]$--module by an element $w_{\lambda,r}$ with relations:\begin{gather*} \lie n^+[t]w_{\lambda,r}=0,\qquad (x_{i}^-)^{\lambda(h_{i})+1}w_{\lambda,r}=0,\\
 (h\otimes t^s)w_{\lambda,r}= \delta_{s,0}\lambda(h)w_{\lambda,r},
 \end{gather*} where $i\in I$, $h\in\lie h$ and $s\in\bz_+$.
 The following proposition summarizes the properties of $\Delta(\lambda,r)$ which are necessary for this paper (see for example \cite{BCM}).

  \begin{prop} Let $(\lambda,r)\in P^+\times\bz$.
 \begin{enumerit}
 \item[(i)] The module $\Delta(\lambda,r)$ is indecomposable and finite--dimensional and hence an object of $\cal I_{\bdd}$.
 \item[(ii)] $\dim \Delta(\lambda,r)_\lambda =\dim\Delta(\lambda,r)[r]_\lambda=1$,
 \item[(iii)] $\wt\Delta(\lambda,r)\subset\conv W\lambda$,
 \item[(iv)] The module $V(\lambda,r)$ is the unique irreducible quotient of $\Delta(\lambda,r)$.
 \item[(v)]   $\{\ch_{\gr}\Delta(\lambda,r):(\lambda,r)\in P^+\times\bz\}$  is a linearly independent subset of $\bz[P][u,u^{-1}]$.
\end{enumerit}\hfill\qedsymbol
\end{prop}
We denote by $[\Delta(\lambda,r): V(\mu,s)]$ the multiplicity of $V(\mu,s)$ in a Jordan--Holder series of $\Delta(\lambda,r)$.

\subsection{}\label{nabla} We now define the modules  $\nabla(\lambda,r)$. These modules are usually defined to be the dual of the modules $\Delta(\lambda,r)$, but in our situation the resulting modules would be too small. The correct definition is to take $\nabla(\lambda,r)$ to be the dual of the global Weyl modules $W(\lambda,r)$.  Here  $W(\lambda,r)$  is generated as a $\lie g[t]$--module by an element $w_{\lambda,r}$ with relations:\begin{gather*} \lie n^+[t]w_{\lambda,r}=0,\qquad (x_{i}^-)^{\lambda(h_{i})+1}w_{\lambda,r}=0,\\
 hw_{\lambda,r}=\lambda(h)w_{\lambda,r},
 \end{gather*} where $i\in I$ and $h\in\lie h$. Clearly the module $\Delta(\lambda,r)$ is a quotient of $W(\lambda,r)$ and moreover $V(\lambda,r)$ is the unique irreducible quotient of $W(\lambda,r)$.   It is known (see \cite{CFK} or \cite{CPweyl} )  that $W(0,r)\cong\bc$ and that if $\lambda\ne 0$, the modules $W(\lambda,r)$ are infinite-dimensional and satisfy $$\wt W(\lambda,r)\subset\conv W\lambda.$$  It follows that  if we set  $$\nabla(\lambda,r)= W(-w_0\lambda,-r)^*,$$ then $\nabla(\lambda, r)\in\Ob\cal I_{\bdd}$ and $\soc\nabla(\lambda,r)\cong V(\lambda, r)$. The following proposition summarizes the main results on $\nabla(\lambda,r)$ that are needed for this paper.  These follow from corresponding results on $W(\lambda, r)$, see for example \cite{BCM}.
 \begin{prop} Let $(\lambda,r)\in P^+\times\bz$.
\begin{enumerit}
\item[(i)] The module  $\nabla(\lambda, r)$ is an indecomposable object of $\cal I_{\bdd}$.
\item[(ii)] $\dim\nabla(\lambda,r)[r]_\lambda=1$, and $\dim\nabla(\lambda,r)[s]_\lambda\ne 0\iff s\le r$,
\item[(iii)] $\wt\nabla(\lambda,r)\subset\conv W\lambda$,
\item[(iv)] Any submodule of $\nabla(\lambda,r)$ contains $\nabla(\lambda,r)[r]_\lambda$ and the socle of $\nabla(\lambda,r)$ is $V(\lambda,r)$.
     \item[(v)]  $\{\ch_{\gr}\nabla(\lambda,r):(\lambda,r)\in P^+\times\bz\}$ is a linearly independent subset of $\bz[P][[u,u^{-1}]]$.\end{enumerit}\hfill\qedsymbol
 \end{prop}

\subsection{}\label{tilting} \begin{defn}
  We say that $M\in\Ob\cal I$ admits a $\Delta$ (resp. $\nabla$)--filtration if there exists an increasing family of submodules $$0\subset M_1\subset M_2\subset \cdots,\qquad \ M=\bigcup_k M_k,$$  such that\begin{gather*} M_k/M_{k-1}\cong \bigoplus_{(\lambda,r)\in P^+\times\bz}\Delta(\lambda,r)^{m_k(\lambda,r)}, \qquad ({\rm resp.},  \ M_k/M_{k-1}\cong \bigoplus_{(\lambda,r)\in P^+\times\bz}\nabla(\lambda,r)^{m_k(\lambda,r)}), \end{gather*}
 for some choice of $m_k(\lambda,r)\in \bz_+$. We say that $M$ is tilting if $M$ has both a $\Delta$ and a $\nabla$--filtration.
\hfill\qedsymbol
\end{defn}

Since $\dim M[r]_\lambda<\infty$ for all $(\lambda,r)\in P^+\times\bz$, we see that if $M$ has a  $\Delta$--filtration (resp.$\nabla$--filtration) $M_k\subset M_{k+1}$,  then $m_k(\lambda,r)=0$ for all but finitely many $k$.
Since $$\ch_{\gr}M= \sum_{k\ge 0}\ch_{\gr}M_k/M_{k-1}=
\sum_{(\lambda,r)\in\bz}\left(\sum_{k\ge 0}m_k(\lambda,r)\right)\ch_{\gr}\Delta(\lambda,r),$$ (where we understand that $M_{-1}=0$) it follows from
Proposition \ref{delta} that the filtration multiplicity  $$[M: \Delta(\lambda,r)]=\sum_{k\ge 0}m_k(\lambda,r),
$$ is well -defined and  independent of the choice of the filtration. An analogous statement holds for modules admitting a $\nabla$--filtration.

\subsection{} The main goal of this paper is to understand tilting modules in $\cal I_{\bdd}$. In the case of algebraic groups (see \cite{Donkin}, \cite{Mathieu})  a crucial necessary result is to give a cohomological characterization of  modules admitting a $\nabla$--filtration. The analogous result in our situation is  to prove the following statement:

\vskip 12pt

 An object $M$ of $\cal I_{\bdd}$ admits a $\nabla$--filtration iff $\Ext^1_{\cal I}((\Delta(\lambda,r), M)=0$ for all $(\lambda,r)\in P^+\times\bz$.

\vskip 12pt

It is not hard to see that the forward implication is true. The converse statement however requires one to prove  that any object of $\cal I_{\bdd}$ be embedded in a module which admits a $\nabla$--filtration. At this point we can only prove the result for $\lie{sl}_{n+1}$ and we explain the reason for these limitations  in the next  section. Summarizing, the first main result that we shall prove in this paper is:
\begin{prop}\label{extnablaconn}  Let $M\in\Ob\cal I_{\bdd}$.
 \begin{enumerit}
 \item[(i)] If  $M$ admits a $\nabla$--filtration, then  for all $(\lambda,r)\in P^+\times \bz$, we have $$\Ext^1_{\cal I}(\Delta(\lambda,r), M)=0.$$
\item[(ii)] Let  $\lie g$ be of type $A_n$, and assume that   $M\in\cal I_{\bdd}$  satisfies
 $\Ext^1_{\cal I}(\Delta(\lambda,r), M)=0$ for all $(\lambda,r)\in P^+\times\bz$. Then $M$ admits a $\nabla$--filtration. \end{enumerit}\end{prop}

 \subsection{} The second main result that we shall prove in this paper is the following:
 \begin{thm}\label{mainthm}
 \begin{enumerit}
\item[(i)] Given $(\lambda, r)\in P^+\times\bz$,  there exists an indecomposable module $T(\lambda,r)\in\Ob\cal I_{\bdd}$  which admits a $\Delta$--filtration and satisfies
\begin{gather*} \Ext^1_{\cal I}(\Delta(\mu,s), T(\lambda,r))=0,\ \ (\mu,s)\in P^+\times\bz,\\ T(\lambda,r)[r]_\lambda=1,\ \ \wt T(\lambda,r)\subset\conv W\lambda,\end{gather*} and $T(\lambda,r)\cong T(\mu,s)$ iff $(\lambda,r)=(\mu,s)$.
    \item[(ii)] If  $\lie g$ is of type $\lie{sl}_{n+1}$, then $T(\lambda,r)$ is tilting. Moreover any indecomposable tilting module in $\cal I_{\bdd}$ is isomorphic to $T(\lambda,r)$ for some $(\lambda,r)\in P^+\times\bz$. Finally any tilting module in $\cal I_{\bdd}$ is isomorphic to a direct sum of indecomposable tilting modules.

    \end{enumerit}
    \end{thm}
\section{The canonical filtration and proof of Proposition \ref{extnablaconn}}
In this section we show that one can define in a canonical way a filtration on any object of $\cal I_{\bdd}$ such that the successive quotients embed into a direct sum of modules $\nabla(\mu,s)$, $(\mu,s)\in P^+\times\bz$. To do this we need to understand the projective and injective objects of $\cal I$ although these are not objects of $\cal I_{\bdd}$. Using the canonical filtration we get an upper bound for the character of any object of $\cal I_{\bdd}$. We then  use this bound along with the  BGG--reciprocity result proved in  \cite{BBCKL} and \cite{BCM} to establish  Proposition \ref{extnablaconn}.

\subsection{} The category $\cal I$ contains the projective cover  and the  injective envelope of a simple object. For $(\lambda,r)\in P^+\times\bz$, set
$$ P(\lambda,r)=\bu(\lie g[t])\otimes_{\bu(\lie g)} V(\lambda,r),\qquad \ I(\lambda,r) = P(-w_0\lambda,-r)^*.$$ Note that \begin{gather*} P(\lambda,r)[r]\ \cong_{\lie g} V(\lambda)\  \cong_{\lie g} I(\lambda,r)[r] ,\\ \\ P(\lambda, r)[s]\ =0\  = I(\lambda, -r)[-s] , \ \ \ s<r.\end{gather*} Clearly $P(\lambda,r)$ is generated by the element $p_{\lambda,r}=1\otimes v_{\lambda}$ with defining relations: $$\lie n^+p_{\lambda,r}=0,\ \ h p_{\lambda,r}=\lambda(h)p_{\lambda,r},\ \ (x_i^-)^{\lambda(h_i)+1} p_{\lambda,r}=0.$$
The following  was proved in \cite[Proposition 2.1]{CG1}.
\begin{prop}\label{pdefrel} For $(\lambda,r)\in P^+\times\bz$,  the object $P(\lambda,r)$ is the  projective cover in $\mathcal I$ of  $V(\lambda, r)$. Analogously, the object $I(\lambda,r)$ is the injective envelope of $V(\lambda,r)$ in $\cal I$.  \hfill\qedsymbol
\end{prop}
 Notice that $P(\lambda,r)_\mu\ne 0$ for infinitely many $\mu\ge\lambda$ and hence $P(\lambda,r)$ (and also $I(\lambda,r)$) is not an object of $\cal I_{\bdd}$. However, we shall introduce quotients (resp. submodules) of these objects which do lie in $\cal I_{\bdd}$.

\subsection{} The object  $W(\lambda,r)$ defined in Section \ref{nabla}  is the
the maximal quotient (in $\cal I$)  of $P(\lambda,r)$ such that $$\wt W(\lambda,r)\subset\lambda- Q^+,$$ or equivalently the maximal quotient whose weights are contained in $\conv W\lambda$. Similarly, $\nabla(\lambda,r)$ is the maximal submodule of $I(\lambda,r)$ whose weights are in the $\conv W\lambda$. The following is now trivially proved.
\begin{lem}\label{zeroextglob} For $\lambda,\mu\in P^+$ with $\lambda\nless \mu$,  we have $$\Ext^1_{\cal I}(W(\lambda,r), W(\mu,s))=0 = \Ext^1_{\cal I}(\nabla(\mu,r), \nabla(\lambda,s)),\ \ {\rm{for\ all}}\ \ r,s\in\bz.$$\hfill\qedsymbol
\end{lem}

\subsection{} At this stage it is worth making the following remark. Define a partial order $\preceq$ on $P^+\times \bz$ by: $(\lambda,r)\preceq (\mu,s)$ if either $\lambda<\mu$ or $\lambda=\mu$ and $r\le s$. Then it is not hard to see that, $\Delta(\lambda,r)$ is the maximal quotient of $P(\lambda,r)$ such that $$\Delta(\lambda,r)[s]_\mu\ne 0\implies (\mu,s)\preceq (\lambda,r).$$ On the other hand, $\nabla(\lambda,r)$ is the maximal submodule of $I(\lambda,r)$ satisfying, $$\nabla(\lambda,r)[s]_\mu\ne 0\implies (\mu,s)\preceq (\lambda,r),$$ and hence our choices are consistent with the ones usually made in the literature.

\subsection{}  Given $\Gamma\subset P^+$, let $\cal I(\Gamma)$ be the full subcategory of $\cal I$ consisting of objects $M$ such that
 $$\wt M\subset\bigcup_{\lambda\in \Gamma}\conv W\lambda.$$  The category $\cal I_{\bdd}(\Gamma)$ is defined similarly. Given $M\in\cal I$, let $M_\Gamma$ be the maximal submodule of $M$ that lies in $\cal I(\Gamma)$.
 We shall say that a subset $\Gamma$ of $P^+$ is closed with respect to $\le$ if $\lambda\in\Gamma$ and $\mu\le \lambda$ implies $\mu\in\Gamma$.

 \begin{prop}\label{embed} Let $\Gamma\subset P^+$ be closed with respect to $\le$.  Then $I(\lambda,r)_\Gamma$ is an injective object of $\cal I_{\bdd}(\Gamma)$ for $(\lambda,r)\in  \Gamma\times\bz$. Moreover, if  $M\in\Ob\cal I_{\bdd}(\Gamma)$ for some finite closed subset $\Gamma\subset P^+$,  there exists an injective morphism $$M\hookrightarrow\bigoplus_{(\lambda,r)\in\Gamma\times\bz} I(\lambda,r)_\Gamma^{\oplus m(\lambda,r)},\ \ m(\lambda,r)=\dim\Hom_{\cal I}(V(\lambda,r), M).$$
 \end{prop}
 \begin{pf} The first statement is immediate from the fact that
 $I(\lambda,r)$ is injective in $\cal I$ and the observation that if $\pi:M\to N$ is a morphism of objects in $\cal I$ and $M\in\Ob\cal I_{\bdd}(\Gamma)$, then $\pi(M)\in\Ob{\cal I}_{\bdd}(\Gamma)$. For the second statement let $(\lambda,r)\in\Gamma\times\bz$ and $M\in\Ob{\cal I}_{\bdd}(\Gamma)$.
  Corresponding to any non--zero morphism $\varphi: V(\lambda,r)\to M$,  we have a morphism $\tilde\varphi: M\to I(\lambda,r)$ whose image is clearly in $I(\lambda,r)_\Gamma$. If $m(\lambda,r)>0$ it follows that  by fixing a basis for $\Hom_{\cal I}(V(\lambda,r), M)$ we have a morphism $$\varphi_{\lambda,r}: M\to\left( I(\lambda,r)_\Gamma\right)^{\oplus m(\lambda,r)}.$$  Notice that $$\varphi_{\lambda,r}M[s]\ne 0\implies s\le r.$$
  Since $M[\ell]=0$ for all $\ell>>0$,  it follows that we have  a well-defined map $$\Phi: M\to\bigoplus_{(\lambda,r)\in \Gamma\times\bz} \left(I(\lambda,r)_\Gamma\right)^{m(\lambda,r)},\qquad  \ m\to\{\varphi_{\lambda,r}(m)\}_{(\lambda,r)\in\Gamma\times\bz}.$$ It remains to prove that $\Phi$ is injective. If $\ker\Phi\ne 0$ then we have $\soc\ker\Phi\ne 0$ by Proposition \ref{simpleclass}. On the other hand, $\soc\ker\Phi\subset\soc M$ and the restriction of $\Phi$ to $\soc M$ is injective by design. The proof is complete.
 \end{pf}

\subsection{}\label{enump} {\em From now on we fix an enumeration $\lambda_0,\lambda_1,\cdots,\lambda_k,\cdots$ of $P^+$ satisfying: $$\lambda_r-\lambda_s\in Q^+\implies r\ge s.$$} Given $M\in\cal I_{\bdd}$, define  $k(M)\in\bz_+$ to  be minimal such that $$\wt M\subset\bigcup_{s=0}^{k(M)}\conv  W\lambda_s.$$ For $0\le s\le k(M)$, let $M_s$ be the maximal submodule of $M$ whose weights lie in the union of the sets $\{\conv W\lambda_r: r\le s\}$. Clearly $$M_s\subset M_{s+1},\quad \ M=\bigcup_{s=0}^{k(M)}M_s,\quad  \  \Hom_{\lie g}(V(\lambda, r), M_{s+1}/M_s)\ne 0\implies \lambda=\lambda_{s+1}.$$ We call the filtration $M_0\subset M_1\subset\cdots\subset M_{k(M)}=M$ the canonical filtration of $M$.
  It follows from Proposition \ref{embed}, and the observation that $\Hom_{\lie g}(V(\lambda, r), M_{s+1}/M_s)\ne 0\implies \lambda=\lambda_{s+1}$,  that $M_{s+1}/M_s$ embeds into a direct sum of modules of the form $\nabla(\lambda_{s+1}, r)$, $r\in\bz$, and in fact we get 
\begin{gather}\label{cheq} \ch_{\gr}M=\sum_{s\ge 0}\ch_{\gr}M_{s}/M_{s-1}\le \sum_{s\ge 0}\sum_{r\in\bz}\dim\Hom_{\cal I}(V(\lambda_{s},r), M_s/M_{s-1})\ch\nabla(\lambda_{s},r). \end{gather} 
We claim that this is equivalent to, 
\begin{gather} \label{chineq} \ch_{\gr}M\le \sum_{s\ge 0}\sum_{r\in\bz}\dim\Hom_{\cal I}(\Delta(\lambda_{s},r), M)\ch\nabla(\lambda_{s},r).\end{gather} 
For the claim, observe that any non--zero map $\varphi: \Delta(\lambda_s,r)\to M$ has its image in $M_s$. Moreover $\varphi$ maps the unique maximal, proper submodule of $\Delta(\lambda_s,r)$ to $M_{s-1}$ and hence induces a non--zero map from $V(\lambda_s,r)\to M_s/M_{s-1}$, which proves that there is an injective map from $\Hom_{\cal I}(\Delta(\lambda_{s},r), M)$ to $\Hom_{\cal I}(V(\lambda_{s},r), M_s/M_{s-1})$, and hence that
$$\dim\Hom_{\cal I}(\Delta(\lambda_{s},r), M)\le \dim\Hom_{\cal I}(V(\lambda_{s},r), M_s/M_{s-1}).$$ 
For the reverse inequality, suppose that we have a non--zero map $\psi: V(\lambda_s,r)\to M_s/M_{s-1}$ and choose $m\in  M_s[r]_{\lambda_s}$ such that  $\psi(v_{\lambda_s,r})= \bar m$ where $\bar m$ is the image of $m$ in $M_s/M_{s-1}$. Since $\wt M_s\subset \conv W\lambda_s$ it follows that 
$$\lie n^+[t] m= 0.$$
On the other hand since  $(\lie h\otimes t\bc[t]) \bar m =0,$  we must have that 
$$(\lie h\otimes t\bc[t]) m \in (M_{s-1})_{\lambda_s}=0.$$
Hence there exists a non--zero map from $\Delta(\lambda_s,r)\to M_s$ which proves there is an injective map from $\Hom_{\cal I}(V(\lambda_{s},r), M_s/M_{s-1})$ to $\Hom_{\cal I}(\Delta(\lambda_{s},r), M)$, proving the claim about dimensions.
Finally, note that equality holds in \eqref{chineq} iff the canonical filtration is a $\nabla$--filtration.

 \subsection{} The following result was proved in \cite{BCM} when $\lie g$ is of type $\lie{sl}_2$ and in \cite{BBCKL} when $\lie g$ is of type $\lie{sl}_{n+1}$. More precisely the dual of the following result was proved in these papers, i.e. it was shown that  the projective objects had a canonical decreasing filtration with successive quotients being the global Weyl modules $W(\mu,s)$. It is conjectured in \cite{BCM} that the result is true in general.
 \begin{thm} \label{bgg} Assume that $\lie g$ is of type $\lie{sl}_{n+1}$. Let $\Gamma$ be a finite subset of $P^+$. For all $(\lambda,r)\in\Gamma\times\bz$ the canonical filtration of  $I(\lambda,r)_\Gamma$ is a $\nabla$--filtration. Moreover for all $(\mu,s)\in P^+\times\bz$, we have
 \begin{equation}[I(\lambda,r)_\Gamma:  \nabla(\mu,s)]=[\Delta(\mu,s): V(\lambda,r)]=\dim\Hom_{\cal I}(\Delta(\mu,s), I(\lambda,r)_\Gamma).\end{equation}\hfill\qedsymbol
 \end{thm}

\subsection{} We note the following consequence Proposition \ref{embed} and Theorem \ref{bgg}.
\begin{prop}\label{embedM}  Assume that $\lie g$ is of type $\lie{sl}_{n+1}$ and let $M\in\Ob{\cal I}_{\bdd}$. Then $M$ embeds into  an object $I(M)$ of $\cal I_{\bdd}$ which  admits a $\nabla$--filtration. \hfill\qedsymbol
\end{prop}

\subsection{}\label{cheq1} To prove (ii) of Proposition \ref{extnablaconn}, suppose that $M\in\Ob\cal I_{\bdd}$ satisfies 
$$\Ext^1_{\cal I}(\Delta(\lambda,r), M)=0,\ \ (\lambda, r)\in P^+\times\bz.$$ 
Assume also that we have an embedding $$0\to M\to I(M)\to Q\to 0,$$ where $I(M)\in\Ob\cal I_{\bdd}$ has a $\nabla$--filtration, in which case   $Q\in\Ob\cal I_{\bdd}$. (In particular Proposition \ref{embedM} shows that we can do this when $\lie g$ is of type $\lie{sl}_{n+1}$). Applying $\Hom_{\cal I}(\Delta(\lambda,r),-)$ to the short exact sequence shows that for all $\lambda$
$$\dim\Hom_{\cal I}(\Delta(\lambda,r),I(M))-\dim\Hom_{\cal I}(\Delta(\lambda,r),Q)=\dim\Hom_{\cal I}(\Delta(\lambda,r),M).$$ 
Since 
$$\ch_{\gr} M=\ch_{\gr} I(M) -\ch_{\gr} Q$$
and \eqref{chineq} implies that this is equal to
\begin{gather}{\sum_{s\ge 0} \sum_{r\in \bz} (\dim\Hom_{\cal I}(\Delta(\lambda_s,r), I(M)) -\dim\Hom_{\cal I}(\Delta(\lambda_s,r), Q))\ch_{\gr} \nabla(\lambda_s,r) }\end{gather}
\begin{gather}{= \sum_{s\ge 0} \sum_{r\in \bz}\dim\Hom_{\cal I}(\Delta(\lambda_s,r), M)\ch_{\gr} \nabla(\lambda_s,r)}\end{gather}
part (ii) now follows because \eqref{chineq} is an equality.

\subsection{} We need one more standard result (whose proof we include for convenience) to prove part (i) of Proposition \ref{extnablaconn}. 
\begin{lem}\label{reorder}  Suppose that  $M\in\cal I_{\bdd}$ has a (possibly infinite) $\nabla$--filtration. Then $M$ admits a $\nabla$--filtration 
$$0\subset M_1\subset M_2\subset \cdots\subset M_k=M,\qquad M_s/M_{s-1}\cong\bigoplus_{r\in\bz}\nabla(\lambda_s,r)^{\oplus[M:\nabla(\lambda_s,r)]}$$ 
where we recall that $[M:\nabla(\lambda_s,r)]<\infty$ for all $s$ and $r$.  In particular there exists $(\mu,s)\in P^+\times \bz$ with $\mu\in P^+$ maximal such that $M_\mu\ne 0$  and a surjective map $M\to\nabla(\mu,s)$ such that the kernel of this map also admits a $\nabla$--filtration.
 \end{lem}

\begin{pf} Let $N_\ell\subset N_{\ell+1}$ be a $\nabla$--filtration of $M$ and assume that $\lambda\in P^+$ is minimal such that $[M:\nabla(\lambda,r)]\ne  0$ for some $r\in\bz$.  Using Lemma \ref{zeroextglob} and an induction on $\ell$, we see that $$\Ext^1(\nabla(\lambda,r), N_\ell)=0,\ \ \ell\ge 1,\ \ r\ge \bz.$$ 
This implies that for each $\ell$, we have $\tilde N_\ell\subset N_\ell$ such that 
\begin{gather*}\tilde N_\ell\cap N_{\ell-1}=0,\ \quad \tilde N_\ell\cong\bigoplus_r\nabla(\lambda,r)^{\oplus m_\ell(\lambda,r)},\ \quad \frac{N_\ell}{N_{\ell-1}\oplus\tilde N_\ell}\cong\bigoplus_{(\mu,s):\mu\ne\lambda} \nabla(\mu,s)^{\oplus m_\ell(\mu,s)}.\end{gather*}
Define a  filtration $M_\ell\subset M_{\ell+1}$, $\ell\ge 1$ of $M$ by,
$$M_\ell= N_{\ell-1}\bigoplus_{s>\ell}\tilde N_s,$$ where we recall that $N_0=0$. Then 
$$\frac{M_\ell}{M_{\ell-1}}\cong \frac{N_\ell}{N_{\ell-1}\oplus\tilde N_\ell}\cong\bigoplus_{(\mu,s):\mu\ne\lambda} \nabla(\mu,s)^{\oplus m_\ell(\mu,s)}.$$ 
Because $M$ admits only finitely many dominant integral weights, an iteration of this argument completes the proof.
\end{pf}
\subsection{} The following Lemma establishes Proposition \ref{extnablaconn}(i).
\begin{lem}\label{extdelnab} We have $\Ext^1_{\cal I}(\Delta(\lambda, r), \nabla(\mu,s))=0$ for all $(\lambda,r), (\mu,s)\in P^+\times\bz$. In particular if $N\in\Ob{\cal I}_{\bdd}$ has a  $\nabla$--filtration then $ \Ext^1_{\cal I}(\Delta(\lambda,r), N)=0$.
\end{lem}
\begin{pf} The proof is standard. Thus, suppose  that we have a short exact sequence $$0\to \nabla(\mu,s)\stackrel{\iota}\to M\stackrel{\tau} \to\Delta(\lambda,r)\to 0.$$  Then $M_\lambda\ne 0$ and if $\mu\ngeq\lambda$ we have  $$(\lie n^+[t]) M_\lambda = 0 = (\lie h\otimes t\bc[t]) M_\lambda. $$ It follows from the defining relations of $\Delta(\lambda,r)$  that if $m\in M[r]_\lambda$ is such that $\tau(m)=w_{\lambda,r}$, then $\bu(\lie g[t]) m$ is a quotient of $\Delta(\lambda,r)$ via the map $w_{\lambda,r}\to m$ and hence the sequence splits. If $\mu\ge \lambda$, then  by taking duals we have  a short exact sequence $$0\to \Delta(\lambda,r)^*\stackrel{\tau^*}\to M^*\stackrel{\iota^*}\to W(-w_0\mu, -s)\to 0.$$  Since $-w_0\mu\ge -w_0\lambda$ we have  $\lie n^+[t] M^*_{-w_0\mu}=0$ and  using the defining relations of $W(-w_0\mu, -s)$ we see that $\iota^*$ splits.

 Suppose that $N\in\Ob\cal I_{\bdd}$ admits a $\nabla$--filtration and let $p\in\bz$ be such that $N[s]=0$ if $s>p$. It follows from Lemma \ref{reorder} that there exists $k\in\bz_+$ and a filtration $0\subset N_0\subset N_1\subset\cdots\subset N_k=N$ such that $$N_s/N_{s-1}\cong\bigoplus_{\ell\le p}\nabla(\lambda_s,\ell)^{m(\lambda_s,\ell)},$$ for some $m(\lambda_s,\ell)\in\bz_+$.
Since $$\Ext^1_{\cal I}(\Delta(\lambda,r), N_s/N_{s-1}) \hookrightarrow\prod_{s\le p}(\Ext^1_{\cal I}(\Delta(\lambda,r), \nabla(\lambda_s,\ell)^{\oplus m(\lambda_s,\ell)})$$ it follows that $\Ext^1(\Delta(\lambda,r), N_s/N_{s-1})=0$. An obvious induction on $s$ proves the Lemma.
\end{pf}

\subsection{} \begin{prop}\label{cangoodsl} Suppose that $\lie g$ is of type $\lie{sl}_{n+1}$. An object $M$ of $\cal I_{\bdd}$ has a $\nabla$--filtration iff the canonical filtration of $M$ is a $\nabla$--filtration.  \hfill\qedsymbol\end{prop}
\begin{pf} Suppose that $M$ has a $\nabla$--filtration. Then we have proved in Section \ref{cheq1} that $$\ch_{\gr}M=\sum_{(\lambda,r)\in P^+\times\bz}\dim\Hom_{\cal I}(\Delta(\lambda,r), M)\ch_{\gr}\nabla(\lambda,r).$$ Hence equality must hold in \eqref{chineq} which was written for the canonical filtration. This proves that the canonical filtration is a $\nabla$--filtration. The converse is obvious.

\end{pf}

\section{Modules with $\Delta$--filtrations}
In our situation the fact that the dual of a $\Delta$--module is not a $\nabla$--module means that we have to also study properties of modules admitting a $\Delta$--filtration. We also need some results on the vanishing of  $\Ext^1_{\cal I}(\Delta(\lambda,r), \Delta(\mu,s))$ which will be used to construct the tilting modules in the next section.

\subsection{} Consider the projection map $\bop\bor: \bu(\lie g[t])\to\bu(\lie h[t])\to 0$ corresponding to the vector space decomposition, $$\bu(\lie g[t])= \bu(\lie h[t])\bigoplus\left(\lie n^-[t]\bu(\lie g[t]) +\bu(\lie g[t])\lie n^+[t]\right).$$ For $i\in I$, define elements $P_{i,s}\in\bu(\lie h[t])$ recursively, by $$P_{i,0}=1,\qquad P_{i,s}=-\frac 1s\sum_{r=1}^s(h_i\otimes t^r)P_{i,s-r}.$$ The following was proved in \cite{Garland} (see \cite{CPweyl}) for the current formulation:
\begin{lem}\label{gar} For $i\in I$ and $s\ge 1$, we have,
$$\bop\bor((x_i^+\otimes t)^{s}(x_i^-)^{s})=(-1)^s (s!)^2 P_{i,s}.$$\hfill\qedsymbol
\end{lem}

 \subsection{} \begin{prop}\label{lamnabext}
  \begin{enumerit}
  \item[(i)] Let $(\lambda,r)\in P^+\times\bz$ and assume  that $N\in\Ob{\cal I}$ satisfies,
\begin{equation}\label{cond} N[s]_\lambda =0 \ \ {\rm {if}}\ \ r\le s\le r+1+\sum_{i=1 }^n\lambda(h_i).\end{equation} Then, $$\Ext^1_{\cal I}(\Delta(\lambda,r), N)=0.$$
\item[(ii)] If $\lambda, \mu\in P^+$ and $\mu\ngeq\lambda$, we have  $$\Ext^1_{\cal I}(\Delta(\lambda,r),\Delta(\mu,\ell)) =0,\ \ {\rm{for\ all}}\ \ r, \ell\in\bz,$$ and $$\Ext^1_{\cal I}(\Delta(\lambda,r),\Delta(\lambda, r)) =0,\ \ {\rm{for\ all}}\ \ r\in\bz,$$
    \item[(iii)] Given $\lambda,\mu\in P^+$ there exists $d(\lambda,\mu)\in\bz_+$ such that $$\Ext^1_{\cal I}(\Delta(\lambda,r),\Delta(\mu,s))\ne 0\implies | r-s|\le d(\lambda,\mu).$$
    \end{enumerit}
    \end{prop}
\begin{pf}
  Consider a short exact sequence, $$0\to N\stackrel{\iota}\to M\stackrel{\tau} \to\Delta(\lambda,r)\to 0.$$ Choose $m\in M[r]_\lambda^{\lie n^+}$ such that $\tau(m)=w_{\lambda,r}$. Then $\tau((h_i\otimes t^s)m)=0$ for all $s>0$ or equivalently $(h\otimes t^s)m\in N$.  Using equation \eqref{cond} we get $$ (h_i\otimes t^s)m=0,\ \qquad  0< s\le 1+\sum_{i=1}^n\lambda(h_{i}).$$  Taking $s=1$ gives $$2(x_{\alpha_i}^+\otimes t)m= [h_i\otimes t,x_{\alpha_i}]=0,$$ and repeating we find that for all all $i\in I$ and $k\in\bz_+$ we have $(x_{\alpha_i}^+\otimes t^k)m=0.$ Applying  Lemma \ref{gar}  we have  $$(x^+_{\alpha_i}\otimes t)^{s}(x^-_{\alpha_i})^{s}m= P_{i,s}m=0,\ \qquad s>\lambda(h_i).$$ Since $P_{i,s}$ is a polynomial in $h_i\otimes t^k$, $1\le k\le s$, it follows by an obvious induction that $(h_i\otimes t^s)m=0$ for all $i$ and $s$. Hence we have proved that $m$ satisfies the defining relations of $\Delta(\lambda,r)$ which means that $\tau$ splits.

   The proof of (ii) is similar and easier and we omit the details. Part (iii) is immediate from part (i) and the fact that $\Delta(\mu,s)$ is finite--dimensional.

  \end{pf}

\begin{cor}\label{deltacrit}  If  $M\in\cal I$ has a $\Delta$--filtration then $\Ext^1_{\cal I}(M,\nabla(\lambda,r))=0$ for all $(\lambda,r)\in P^+\times\bz$.
\end{cor}
\begin{pf} If $M$ has a finite $\Delta$--filtration  then a obvious induction on the length of the filtration gives the result. The proof of the infinite  case is a simple exercise and we omit the details.
\end{pf}

 \iffalse

  \begin{prop} Suppose that $N\in\Ob\cal I$ is such that $\Ext^1(\Delta(\lambda,r), N)=0$ for all $(\lambda,r)\in P^+\times \bz$. If $M$ has a $\Delta$--filtration then $\Ext^1(M, N)=0$. \end{prop}
\begin{pf}
 Consider a short exact sequence $$0\to N\to U\to M\to 0.$$ Suppose that $M_k\subset M_{k+1}$ is a part of the $\Delta$-filtration of $M$ and assume that 
 $$M_{k+1}/M_k\cong\Delta(\mu_k,s_k),\ \ {\rm{for\ some}}\ \ (\mu_k,s_k)\in P^+\times \bz.$$
Let $U_k\subset U$ be the preimage of $M_k$ and note that $U_{k+1}/U_k\cong \Delta(\mu_k,s_k)$.
Consider the short exact sequence $$0\to N\to U_k\to M_k\to 0.$$ This sequence defines an element of $\Ext^1(M_k, N)$. Since $M_k$ has a finite $\Delta$--filtration it follows that $\Ext^1(M_k, N)=0$.  Hence the sequence splits and we have a map $\varphi_k: U_k\to N$. We want to prove that $\varphi_{k+1}:U_{k+1}\to N$ can be chosen to extend $\varphi_k$. For this, applying $\Hom(--,N)$ to 
$$0\to U_k\to U_{k+1}\to\Delta(\mu_k,s_k)\to 0,$$
we get 
$$\Hom(U_{k+1}, N)\to \Hom(U_k, N)\to\Ext^1(\Delta(\mu_k,s_k), N)=0,$$ 
which shows that we can choose $\varphi_{k+1}$ to lift $\varphi_k$. Now defining 
$$\varphi: U\to N,\ \ \varphi(u)=\varphi_k(u), \ \ u\in U_k,$$ 
we have the desired splitting of the original short exact sequence.
\end{pf}
Together with Lemma \ref{extdelnab} we now have,

\fi
\subsection{}
\begin{lem}{\label{homext}} Let $(\lambda,r), (\mu,s)\in P^+\times\bz$.
 We have $$\Hom_{\cal I}(\Delta(\lambda,r), \nabla(\mu,s))\cong\begin{cases}\bc,\ \  (\lambda,r)=(\mu,s),\\ 0\ \ \ {\rm{otherwise}}.\end{cases}$$
   \end{lem}
    \begin{pf} Suppose that $\varphi:\Delta(\lambda,r)\to\nabla(\mu,s)$ is non--zero. Then $\varphi(w_{\lambda,r})\ne 0$ and hence we have $\lambda\le \mu$ and $r\le s$. Moreover since any submodule of $\nabla(\mu,s)$ has non--zero socle it follows that $\nabla(\mu,s)[s]_\mu$ must be in the image of $\varphi$ which shows that $\mu\le\lambda$ and $s\ge r$.

    \end{pf}
   \subsection{}  We end the section with a  final result needed to construct $T(\lambda,r)$. It can be deduced from the fact (proved in \cite{CG1}) that the space of extensions between irreducible objects of $\cal I$ is finite--dimensional, but we include a proof for convenience.
   \begin{prop}\label{extfindim}  For $(\lambda,r), (\mu,s)\in P^+\times\bz$, we have $\dim\Ext^1_{\cal I}(\Delta(\lambda,r), \Delta(\mu,s))<\infty$.
   \end{prop}
   \begin{pf}  Let $\pi: P(\lambda,r)\to \Delta(\lambda,r)\to 0$ be the canonical projection which maps $p_{\lambda,r}$ to $w_{\lambda,r}$. Apply $\Hom(--,\Delta(\mu,s))$ to the short exact sequence $$0\to \ker\pi\to P(\lambda,r)\to\Delta(\lambda,r)\to 0.$$ Since $P(\lambda,r)$ is a projective object of $\cal I$, the result follows if we prove that \begin{equation}\label{homfindim} \dim\Hom_{\cal I}(\ker\pi,\Delta(\mu,s))<\infty.\end{equation} Choose $\ell\in\bz$ such that $\Delta(\mu,s)[p]=0$ for all $p>\ell$, in which case we have an injective map $$\Hom_{\cal I}(\ker\pi,\Delta(\mu,s))\to \dim\Hom_{\cal I}\left(\frac{\ker\pi}{\bigoplus_{p>\ell}\ker\pi[\ell]},\Delta(\mu,s)\right).$$ Since $$\dim\left(\frac{\ker\pi}{\bigoplus_{p>\ell}\ker\pi[\ell]}\right)=\sum_{p=r}^\ell\dim\ker\pi[p]<\infty,$$ equation \ref{homfindim} is proved.
\end{pf}

\section{The  modules $T(\lambda,r)$}
In this section we construct a family of indecomposable modules $T(\lambda,r)$, $\lambda\in P^+$, $r\in\bz$, satisfying: 
$$\dim T(\lambda,r)[r]_\lambda=1,\qquad T(\lambda,r)_\mu\ne 0\implies \mu\le\lambda.$$ 
The construction is similar to the one given in \cite{Mathieu} but there are several difficulties to be overcome in our situation.  By using the operators $\tau_r$, it is enough to construct $T(\lambda,0)$.  We note that our construction makes frequent use of the ordering established in Section 3.5, and so we fix $k\ge 0$ such that $\lambda=\lambda_k$.

\subsection{} We begin by noting the following elementary result.
\begin{lem}\label{elem} Suppose that $M,N\in\Ob\cal I$ are such that $0<\dim\Ext^1_{\cal I}(M,N)<\infty$ and $\Ext^1_{\cal I}(M,M)=0$. Then, there exists $U\in\Ob\cal I$, $d\in\bz_+$ and a non--split short exact sequence 
$$0\to N\to U\to M^{\oplus d}\to 0$$ 
so that $\Ext^1_{\cal I}(M,U)=0$.
\hfill\qedsymbol\end{lem}
%\begin{rem}\label{split} If $\Ext^1(M,N)=0$, then it is convenient to set $\be(M,N)= M$.\end{rem}

%\subsection{} Set $r_k=0$ and for $0<s\le k$, choose $r_s\in\bz$ with $0<r_{k-1}<r_{k-2}<\cdots<r_1$ such that the following two conditions are satisfied:
%\begin{equation}\label{eq1} \Delta(\lambda_s, r_s)[\ell]\ne 0 \implies \Delta(\lambda_p, r_p)[\ell]=0,\ \qquad 0\le p<s\le k.
%\end{equation}
%This  choice can be made since $\Delta(\lambda,r)$ is finite--dimensional and so has only finitely many graded pieces.  In particular, what we are doing here is to choose the integers $r_s$ so that the following holds. Set  $$V=\bigoplus_{s=0}^k\Delta(\lambda_s, r_s). $$ Then $V[p]\ne 0$ implies that there exists an unique $0\le s\le k$ such that $\Delta(\lambda_s,r_s)[p]\ne 0$. Moreover, if $p<p'$ and $V[p']$ is also non--zero, and $s'$ is such that $\Delta(\lambda_{s'},r_{s'})[p']\ne 0$ then we must have $s'<s$.
%Proposition \ref{lamnabext}(iii) gives and

Set $r_k=0$ and for $0\le s<k$ recursively define $r_s \ge r_{s+1}$ by setting 
\[r_s=\text{max}\{r\ge r_{s+1}: \Delta(\lambda_{s+1}, r_{s+1})[r]\ne 0\}.\]  
It is easy to see that this choice implies
\begin{equation}\label{crucial}\Ext^1_{\cal I}(\Delta(\lambda_s, r),\Delta(\lambda_p, \ell))=0\ \ {\rm{for\ all}}\ \ s<p, \ \ r_s<r,\ \ \ell\le r_p.\end{equation}  

Consider the set $\cal S=\{(\lambda_s,r): 0 \le s \le k,  \ r\le r_s\}$ and let  $\eta: \cal S\to\bz_+$ be the enumeration given by 
$$\eta(\lambda_s, \  r_s-\ell)= k-s+ (k+1)\ell,\ \  s\ge 0.$$
Given $s\in\bz_+$, let  $(\mu_s,p_s)\in\cal S$, be the unique element such that 
$$\eta(\mu_s,p_s)=s,\qquad (\mu_0,p_0)=(\lambda, 0).$$ 
%Since $\tau_r(\Delta(\mu,s))=\Delta(\mu, r+s)$, we find that
%\begin{equation}\label{eq2} \Delta(\mu_s, p_s)[p_\ell]_{\mu_\ell}\ne 0\implies s\le \ell.\end{equation}
Proposition \ref{lamnabext}(ii) and \eqref{crucial}, gives
\begin{equation}\label{crucialtoo}\Ext^1_{\cal I}(\Delta(\mu, r), \Delta(\mu_s, p_s))=0,\ \ s\ge 0,\ \ (\mu,r)\notin\cal S.\end{equation} \iffalse 
Moreover, equation \eqref{eq1} gives,
$$\Delta(\mu_s, p_s)[\ell]\ne 0 \implies \Delta(\mu_r,p_r)[\ell]= 0\ \ s<r.$$\fi

\begin{rem}
To illustrate this proceedure, consider the case when $\lie g= \lie{sl_2}$ and $\lambda=4\omega$, which is $\lambda_4$ by the enumeration given in Section 3.5.  By examining the graded characters of the local Weyl modules (see for example \cite{BCM}, Section 5.12) we see that 
\[ r_3=4,\ \ r_2=6,\ \ r_1=7,\ \ r_0=7. \]
The enumeration given by $\eta$ is given by 
\begin{gather*}(\mu_0, p_0)=(4 \omega, 0), \ \ (\mu_1, p_1)=( 3\omega, 4), \ \ ( \mu_2, p_2)=(2 \omega, 6), \\ 
 (\mu_3, p_3)=(\omega, 7), \ \ (\mu_4, p_4)=(0, 7), \ \ (\mu_5, p_5)=(4\omega, -1) \end{gather*}

Note that Proposition \ref{lamnabext} impies that $\Ext^1_{\cal I}(\Delta(i \omega, r), \Delta(4 \omega, s))=0$ if $i$ is odd.
\end{rem}

\subsection{} Proposition \ref{extfindim} implies that if  $N\in\Ob\cal I$ admits a finite $\Delta$--filtration, then $$ \dim\Ext^1_{\cal I}(\Delta(\mu,r), N)<\infty\ \  {\rm{for\ all}}\ (\mu,r)\in P^+\times\bz.$$ 
We now use  Proposition \ref{lamnabext} and  Lemma \ref{elem}  to define finite--dimensional objects $M_s$,  $s\ge 0$, of $\cal I$, recursively as follows.

Set $M_0 =\Delta(\lambda, 0)$. If $\Ext^1_{\cal I}(\Delta(\mu_1,p_1), M_0)=0$, take  $M_0=M_1$. Otherwise,
let $U\in\Ob\cal I_{\bdd}$ be chosen as in Lemma \ref{elem} (with $M= \Delta(\mu_1,p_1)$ and $N=M_0$) and let $M_1$ be the indecomposable summand of $U$ which contains $U[0]_{\lambda}\cong (M_0)_{\lambda}$ and note that $$\Ext^1_{\cal I}(\Delta(\mu_1,p_1), M_1)=0.$$  
Since $U$ is finite--dimensional and $\Delta(\mu_1,p_1)$ is indecomposable, it follows that there exists $d_1\in\bz_+$ and a non--split short exact sequence of the form 
$$0\to M_0\stackrel{\iota_0}\to M_1\stackrel{\tau}\to\Delta(\mu_1,p_1)^{\oplus d_1}\to 0. $$
Clearly, $M_1$ is generated as a $\lie g[t]$--module by the spaces $M_1[0]_{\lambda}$ and $M_1[p_1]_{\mu_1}$,and  
\begin{equation}\label{zero}\Ext^1_{\cal I}(\Delta(\mu_j,p_j),M_1)=0,\ \ j=0,1,\ \   \Ext^1_{\cal I}(\Delta(\mu,r), M_1)=0,\ \ (\mu,r)\notin \cal S.\end{equation}
Repeating this procedure, we can  construct a family $M_s$, $s\ge 0$, of indecomposable finite--dimensional modules and injective morphisms $\iota_s: M_s\to M_{s+1}$ of objects of $\cal I$. Each $M_s$ admits a  finite  $\Delta$--filtration, and satisfies  
$$\dim M_s[0]_{\lambda}=1,\qquad \wt M_s\subset\conv W\lambda,$$ 
and,
$$\Ext^1_{\cal I}(\Delta(\mu_\ell,p_\ell), M_s)=0,\ \ 0\le\ell\le s,\qquad \Ext^1_{\cal I}(\Delta(\mu,r), M_s)=0,\ \ (\mu,r)\notin \cal S.$$

If we choose $\ell_0$ to be maximal such that $\Delta(\mu_k,p_k)[\ell_0]\ne 0$, then $\ell_0$ is also maximal such that $M_s[\ell_0]\ne 0$ for $0\le s\le k$. It follows that  
$$\Delta(\mu_s, p_s)[\ell]\ne 0 \implies \ell\le \ell_0,\ \ s\ge 0,$$ 
and so we have
\begin{equation}\label{bddgrade} M_s[p]=0, \ {\rm{for \ all}}\ \  s\ge 0,\  p>\ell_0\end{equation}
 establishing that $M_s \in \Ob\cal I_{\bdd}$.  Set $\iota_{r,s}= \iota_{s-1}\cdots\iota_r: M_r\to M_{s},\ \ r< s,\ \ \iota_{r,r}=\id. $  Then,   \begin{equation}\label{ms}  M_s[p_\ell]_{\mu_\ell}= \iota_{\ell,s}(M_\ell[p_\ell]_{\mu_\ell}),\ \ s\ge\ell, \end{equation}
and $M_s$ is generated as a $\lie g[t]$--module by the spaces $\{M_s[p_\ell]_{\mu_\ell}: \ell\le s\}$.
Let $ T(\lambda,0) $ be the direct limit of  $\{M_s, \iota_{r,s}: r,s\in \bz_+, r\le s.\}$, and set $T(\lambda,r)=\tau_rT(\lambda,0)$.

It is straightforward to see that the preceding discussion establishes the following.
 \begin{lem} For $(\lambda,r)\in P^+\times \bz$ the module $T(\lambda,r)$ is an object of $\cal I_{\bdd}$. We have $\wt T(\lambda,r) \subset\conv W\lambda$ and $\dim T(\lambda,r)[r]_{\lambda}=1$.  Clearly  $T(\lambda,r)\cong T(\mu, s)$ iff $\lambda=\mu$ and $r=s$.\hfill\qedsymbol
 \end{lem}

\subsection{} Since the maps $\iota_{r,s}$ are injective morphisms it follows that  the canonical morphism $ M_s\to T(\lambda,0)$ is injective and we have an isomorphism of  $M_s$ with a submodule $\tilde M_s$ of $T(\lambda,0)$. Moreover, we have inclusions $\tilde M_s\subset\tilde M_{s+1}$ and 
$$T(\lambda,0)=\bigcup_{s\ge 0}\tilde M_s,\ \qquad\tilde M_{s}/\tilde M_{s-1}\cong M_{s}/M_{s-1},\ \  s\ge 0,$$ proving that $T(\lambda,0)$ has a $\Delta$-filtration.  From now on,  by abuse of notation, we write $M_s$ for $\tilde M_s$. Then, \eqref{ms} gives, 
\begin{equation}\label{ms1} T(\lambda, 0)[p_\ell]_{\mu_\ell} = M_{\ell}[p_\ell]_{\mu_\ell}.\end{equation}
To prove that $T(\lambda,0)$ is indecomposable, suppose that 
$$T(\lambda,r)= U_1\oplus U_2.$$ 
Since $\dim T(\lambda,0)[0]_{\lambda}=1$, we may assume without loss of generality that $T(\lambda,0)[0]_{\lambda}\subset U_1$ and hence $ M_0\subset U_1$. Assume that we have proved by induction that $M_{s-1}\subset U_1$. Since  $M_s$ is generated as a $\lie g[t]$--module by the spaces $\{M_s[p_\ell]_{\mu_\ell}: \ell\le s\}$, it suffices to prove that $M_s[p_s]_{\mu_s}\subset U_1$. By \eqref{ms1}, we have  $U_i[p_s]_{\mu_s}\subset M_s$ and hence  
$$M_s= \left (M_{s-1}+ \bu(\lie g[t])U_1[p_s]_{\mu_s}\right)\bigoplus\bu(\lie g[t])U_2[p_s]_{\mu_s}.$$ 
Since $M_s$ is indecomposable by construction, it follows that $U_2[p_s]_{\mu_s}=0$ and $M_s\subset U_1$ which completes the inductive step.

\subsection{} 
\begin{prop}\label{tiltingis} For all $\lambda\in P^+$ and $(\mu,s)\in P^+\times\bz$, we have
\begin{equation}\label{extzero}\Ext^1_{\cal I}(\Delta(\mu,s)), T(\lambda,0))=0,\end{equation}
\end{prop}

\begin{pf}  Consider a short exact sequence
\begin{equation}\label{ses} 0\to T(\lambda,0)\to U\to \Delta(\mu,s)\to 0.\ \ \end{equation}
If $\mu\nleq\lambda$, an argument identical to the one given in the proof of Lemma \ref{extdelnab} proves that the short exact sequence in \eqref{ses}  must split. If $\mu\le \lambda$ choose $ r>>0$ so that 
\begin{gather} \left(T(\lambda,0)/M_r\right)[\ell]=0,\ \  s\le\ell\le s+1+\sum_{i=1}^n\mu(h_i),\\ \label{const} \Ext^1_{\cal I}(\Delta(\mu,s), M_r)=0.\end{gather} 
We can choose such an $r$  for the following reasons. Since $T(\lambda,0)$ has finite--dimensional graded pieces  there exists $p$ such that for all $r\ge p$ we have $M_r[\ell]= M_p[\ell]$ for all $s\le \ell\le s+1+\sum_{i=1}^n\mu(h_i)$. If $(\mu,s)\notin \cal S$ then \eqref{const} is automatically satisfied. If $(\mu,s)\in \cal S$,  say $\eta((\mu, s))= \tilde s$, then \eqref{const} holds, because of the way  $M_r$ was constructed, if $r> \tilde s$.

Consider the short exact sequence 
$$0\to M_r\to T(\lambda,0)\to T(\lambda,0)/M_r\to 0.$$ 
Applying $\Hom_{\cal I}(\Delta(\mu,s),--)$ to the short exact sequence we get from Proposition \ref{lamnabext}(i) that 
$$\Ext^1_{\cal I}(\Delta(\mu,s), T(\lambda,0)/M_r)=0.$$ 
Using \eqref{const} we see that equation \eqref{extzero} is proved.
\end{pf}

 The following result is now a consequence of  Proposition \ref{extnablaconn}.
\begin{cor} If  $\lie g$ is of type $\lie{sl}_{n+1}$, the objects $T(\lambda,r)$ are tilting.\hfill\qedsymbol \end{cor}

\subsection{} Assume  from now on that Proposition \ref{extnablaconn} is true in which case $T(\lambda,r)$ is tilting. The following result,  which is proved in the rest of the section, completes the proof of Theorem \ref{mainthm}.  It also shows that our construction of tilting modules does not depend on our enumeration of $P^+$.  
\begin{prop}\label{finalprop}  Assume that $T(\lambda,r)$ is a tilting module for all $(\lambda,r)\in P^+\times\bz$. Then any tilting module in $\cal I_{\bdd}$ is isomorphic to direct sum of  modules $T(\lambda,r)$,$ (\lambda,r)\in P^+\times \bz$.
\end{prop}

\subsection{} Let  $T\in\cal I_{\bdd}$ be a fixed tilting module. Using Proposition \ref{extnablaconn} and Corollary \ref{deltacrit}, we have  
\begin{equation}\label{tiltext}  \Ext^1_{\cal I}(T,\nabla(\mu,r))=\Ext^1_{\cal I}(\Delta(\lambda,r), T))=0,\ \ (\lambda,r)\in P^+\times\bz.\end{equation}

\begin{lem}\label{sumtilt} Suppose that $T_1$ is any summand of $T$. Then $T_1$ admits a $\nabla$--filtration and $$\Ext^1_{\cal I}(T_1,\nabla(\lambda,r))=0,$$ for all $(\lambda,r)\in\bz$.
\end{lem}
\begin{pf} Since $\Ext^1$ commutes with finite  direct sums, we get $$\Ext^1_{\cal I}(T_1,\nabla(\lambda,r))=0,\ \ \Ext^1_{\cal I}(\Delta(\lambda,r), T_1))=0,\ \ (\lambda,r)\in P^+\times\bz.$$ Under the assumption that Proposition \ref{extnablaconn}  is true, the second equality implies that $T_1$ has a $\nabla$-filtration and the proof of the Lemma is complete. \end{pf}

\subsection{} The preceding lemma illustrates one of the difficulties we face in our situation. Namely, we cannot directly conclude that $T_1$ has a $\Delta$--filtration from the vanishing $\Ext$--condition by using Proposition \ref{extnablaconn}. However, we can prove,

\begin{prop}\label{tiltext2} Suppose that $N\in\cal I_{\bdd}$ has a $\nabla$--filtration and satisfies 
$$\Ext^1_{\cal I}(N,\nabla(\lambda,r)) =0,\ \ {\rm{for\ all}} \ (\lambda,r)\in P^+\times \bz.$$ 
There exists $(\mu,s)\in P^+\times\bz$ such that $T(\mu,s)$ is a summand of $N$.\end{prop}
\begin{pf}
Since $N$ has a $\nabla$--filtration we can choose $(\mu,s)\in P^+\times\bz$ so that  we have a non--zero surjective map $\varphi: N\to\nabla(\mu,s)\to 0$ and we can  also choose $\pi: T(\mu,s)\to\nabla(\mu,s)\to 0$. We may also assume that $\ker\varphi$ and $\ker\pi$ have $\nabla$--filtrations. Let $v_{\mu,s}$ be a non--zero element of $\nabla(\mu,s)[s]_\mu$ and  choose $m\in N[s]_\mu$ and $u\in T(\mu,s)[s]_\mu$ so that 
$$\varphi(m)=v_{\mu,s}=\pi(u).$$  
Consider the short exact sequences 
$$0\to\ker\varphi\to N\to\nabla(\mu,s)\to 0,$$ 
and 
$$0\to \ker\pi\to T(\mu,s)\to \nabla(\mu,s)\to 0.$$ 
Apply $\Hom_{\cal I}(T(\mu,s),--)$ to the first sequence and $\Hom_{\cal I}(N,--)$ to the second sequence. Since $\ker\varphi$ and $\ker\pi$  admit a  $\nabla$--filtration, equation \eqref{tiltext} gives $\Ext^1_{\cal I}(T(\mu,s), \ker\varphi)\  = 0$. By hypothesis, we also have $\ \Ext^1_{\cal I}(N,\ker\pi) =0 $ and so we have surjective maps \begin{gather*} \Hom_{\cal I}(T(\mu,s), N)\to \Hom_{\cal I}(T(\mu,s),\nabla(\mu,s))\to 0, 
\Hom_{\cal I}(N, T(\mu,s))\to \Hom_{\cal I}(N,\nabla(\mu,s))\to 0.\end{gather*}
Choose $\tilde\varphi\in\Hom_{\cal I}(N, T(\mu,s))$ and $\tilde\pi\in\Hom_{\cal I}( T(\mu,s),N)$ such that
$$\pi.\tilde\varphi=\varphi,\ \ \varphi.\tilde\pi=\pi.$$ 
This gives that 
$$\pi.\tilde\varphi.\tilde\pi=\pi$$ 
Setting $\psi=\tilde\varphi.\tilde\pi$, we see that  $\psi(u)=u$ and  hence $\psi$ is a non--nilpotent endomorphism of $T(\mu,s)$. Moreover, for any $ s$, it follows from \eqref{ms1} that    
$$0\ne \psi(N_s)\subset N_s.$$ 
Since  $N_s$  is indecomposable and finite--dimensional  we can use  Fittings Lemma to conclude that $\psi_s:N_s\to N_s$ is an isomorphism.  It follows that $\psi$ is an isomorphism of $T(\mu,s)$  and hence that $\tilde\pi\psi^{-1}$ is a splitting of $\tilde\varphi:N\to T(\mu,s)$.
\end{pf}

\begin{cor} Any indecomposable tilting module is isomorphic to  $T(\lambda,r)$ for some $(\lambda,r)\in P^+\times\bz$. Further if $T$ is tilting and $(\lambda,r)\in P^+\times\bz$ is such that $T\twoheadrightarrow\nabla(\lambda,r)$ and the kernel admits a $\nabla$--filtration,  then $T(\lambda,r)$ is isomorphic to a direct summand of $T$.\hfill\qedsymbol
\end{cor}
\begin{pf}  Since $T$ is tilting it satisfies \eqref{tiltext} and the corollary follows.
\end{pf}

\subsection{}  Suppose now that $T\in\cal I_{\bdd}$ is a tilting module and let $\lambda\in P^+$ be maximal such that $[T:\nabla(\lambda,r)]\ne 0$ for some $r\in\bz$.  Fix also a decreasing sequence $r_1\ge r_2\ge\cdots $ such that $$[T:\nabla(\lambda,s)]\ne 0\implies s=r_j \ {\rm{for\ some}}\ \ j\ge 1.$$ 
Then we have a surjective map $T\to\nabla(\lambda, r_1)$ and so 
$$T=\iota_1T(\lambda,r_1)\oplus T_1.$$ 
By Lemma \ref{sumtilt}, we see that $T_1$ has a $\nabla$--filtration and that $T_1$ maps onto $\nabla(\lambda, r_2)$ and hence $T(\lambda,r_2)$ is isomorphic to a  summand of $T_1$. Continuing, we find that for $j\ge 1 $, there exists a summand $T_j$ of $T$ with  
$$T=T_j\ \bigoplus_{s=1}^j\iota_s T(\lambda,r_s).$$ 
Let $\pi_j:T\to  \iota_j(T(\lambda,r_j))$ be the canonical projection.  Since $T$ has finite--dimensional graded pieces and $r_j\le r_{j-1}$ are decreasing and the modules $T(\lambda_k,r_j)$ are all graded shifts, it follows that for any $m\in T$ we have $\pi_j(m)=0$ for all but finitely many $j$. Hence we have a surjective map 
$$\pi: T\to \bigoplus_{j\ge 1}\iota_jT(\lambda,r_j)\to 0, \ \ {\rm{and}}\ \ \ker\pi=\bigcap_{j\ge 1} T_j.$$ 
In particular, it follows that
$$T=\ \bigoplus_{j\ge 1}\iota_j T(\lambda,r_j)\oplus \ker\pi.$$ 
Repeat the argument with $\ker\pi$. Since $(\ker\pi)_{\lambda}=0$, the argument stops eventually and Proposition \ref{finalprop} is proved.

\end{document}